\documentclass[12pt,reqno]{amsart}

\usepackage{amsfonts,color,amsthm,amsmath,amssymb}

\usepackage{color}
\def\Box{\vcenter{\vbox{\hrule\hbox{\vrule
     \vbox to 8.8pt{\hbox to 10pt{}\vfill}\vrule}\hrule}}}

\newcommand{\tr}{\text{Tr}}
\def\qed{{\hfill$\square$}}
\def\proof{{\vspace{-0.0cm}\bf Proof: \,}}
\def\N{{\mathbb N}}
\def\Z{{\mathbb Z}}
\def\Q{{\mathbb Q}}

\def\C{{\mathbb C}}
\def\F{{\mathbb F}}
\def\v{\sigma}
\def\mod{{\mathrm{mod\,\,}}}

\def\Tr{{\mathrm{Tr}}}
\def\Norm{{\mathrm{Norm}}}

\def\ord{{\mathrm{ord}}}

\def\lcm{{\mathrm{lcm}}}

\newcommand{\ga}{{\gamma}}

\newtheorem{theorem}{Theorem}[section]

\newtheorem{lemma}[theorem]{Lemma}
\newtheorem{remark}[theorem]{Remark}
\newtheorem{corollary}[theorem]{Corollary}

\newtheorem{proposition}[theorem]{Proposition}
\newtheorem{example}[theorem]{Example}
\newtheorem{conj}[theorem]{Conjecture}
\numberwithin{equation}{section}

\begin{document}
\title[Combinatorial Structures]
{Three-valued Gauss periods, Circulant Weighing Matrices and Association Schemes}

\author[Feng, Momihara and Xiang]{Tao Feng$^*$, Koji Momihara$^{\dagger}$ and Qing Xiang}

\thanks{$^*$Research supported in part by the Fundamental Research Funds for Central Universities of China, the National Natural Science Foundation of China under Grants 11201418 and 11422112, and the Research Fund for Doctoral Programs from the Ministry of Education of China under Grant 20120101120089.}
\thanks{$^{\dagger}$
Research supported by JSPS under Grant-in-Aid for Young Scientists (B) 25800093 and Scientific Research (C) 24540013.}

\address{Department of Mathematics, Zhejiang University, Hangzhou 310027, Zhejiang, P. R. China}
\email{tfeng@zju.edu.cn}

\address{Faculty of Education, Kumamoto University, 2-40-1 Kurokami, Kumamoto 860-8555, Japan} \email{momihara@educ.kumamoto-u.ac.jp}

\address{Department of Mathematical Sciences, University of Delaware, Newark, DE 19716, USA} \email{xiang@math.udel.edu}

\keywords{Association scheme, circulant weighing matrix, cyclotomy, Gauss period, Gauss sum.}

\begin{abstract}
Gauss periods taking exactly two values are closely related to two-weight irreducible cyclic codes and strongly regular Cayley graphs. They have been extensively studied in the work of Schmidt and White and others. In this paper, we consider the question of when Gauss periods take exactly three rational values. We obtain numerical necessary conditions for Gauss periods to take exactly three rational values. We show that in certain cases, the necessary conditions obtained are also sufficient. We give numerous examples where the Gauss periods take exactly three values. Furthermore,  we discuss connections between three-valued Gauss periods and combinatorial structures such as circulant weighing matrices and $3$-class association schemes.
\end{abstract}

\maketitle

\section{Introduction}
Let $\F_q$ be the finite field of order $q$, where $q$ is a power of a prime $p$. Let $\xi_p$ be a complex primitive $p$th root of unity and $\Tr _{q/p}$ be the trace from $\F_q$ to $\Z_p:=\Z/p\Z$. Define
$$\psi\colon\F_q\to \C^{*},\qquad\psi(x)=\xi_p^{\Tr _{q/p}(x)},$$
which is easily seen to be a nontrivial character of $(\F_q, +)$, the additive group
of $\F_q$. Let $\chi\colon\F_q^{*}\to \C^*$ be a multiplicative character of $\F_q$. Define the {\em Gauss sum} by
$$G_q(\chi)=\sum_{a\in {\mathbb F}_q^*}\chi(a)\psi(a).$$
Gauss sums are ubiquitous in number theory and in many areas of combinatorics. Closely related to Gauss sums are the Gauss periods which we define below. As before $q$ is a power of a prime $p$. Let $N>1$ be an integer such that $N|(q-1)$ and $\gamma$ a primitive element of $\F_q$. Then the cosets $C_a^{(N,\F_q)}=\gamma^a \langle \gamma^N\rangle$, $0\leq a\leq N-1$, of $\langle \gamma^N\rangle$ in $\F_q^*$ are called the {\it cyclotomic classes of order $N$} of $\F_q$. We often write $C_a^{(N,q)}$ or simply $C_a$ for $C_a^{(N,\F_q)}$, if there is no confusion.  The corresponding {\it Gauss periods} are defined by
$$\eta_a=\sum_{x\in C_a^{(N,q)}}\psi(x),\; 0\leq a\leq N-1.$$

Even though Gauss sums and Gauss periods were first introduced by Gauss to study cyclotomy  (``circle-splitting"), they have played an important role in the investigations of many combinatorial objects, such as difference sets, irreducible cyclic codes, and strongly regular Cayley graphs, cf. \cite{ehkx, mc, sw, FX, FX1, fmx}. In particular, we note that Gauss sums were used extensively in the work of Baumert and McEliece (\cite{bm}, \cite{mc}) on weights of irreducible cyclic codes. The current paper can be thought as a natural continuation of \cite{sw} in which two-weight irreducible cyclic codes were studied by using Gauss sums. The Gauss periods involved  in \cite{sw} take two distinct rational values as they correspond to the (nonzero) weights of two-weight irreducible cyclic codes. In this paper, we consider Gauss periods which take three distinct rational values, and use them to construct various combinatorial objects such as circulant weighing matrices and association schemes.

A {\it circulant weighing matrix of order $N$} is a square matrix $M$ of the form
\begin{equation}\label{cirweigh}
M=
\begin{pmatrix}
a_0 & a_1 & a_2 & \cdots & a_{N-1}\\
a_{N-1} & a_0 & a_1 & \cdots & a_{N-2}\\
\cdots& \cdots & \cdots & \cdots &\cdots\\
a_1 & a_2 & a_3 & \cdots & a_0\\
\end{pmatrix}
\end{equation}
with $a_i\in \{-1,0,1\}$ for all $i$ and $MM^{\top}=wI$, where $w$ is a positive integer and $I$ is the identity matrix of order $N$. The integer $w$ is called the {\it weight} of the weighing matrix. A circulant weighing matrix of order $N$ and weight $w$ will be denoted by {\bf CW}$(N,w)$.

Let $G$ be an abelian group of order $N$. To facilitate the study of circulant weighing matrices we use the group ring language. The elements of $\C[G]$ are
$$A=\sum_{g\in G}a_gg,$$
with $a_g\in \C$;  for any integer $t$, we write
$$A^{(t)}:=\sum_{g\in G}a_gg^t.$$
For a subset $A$ of $G$, it is customary to identify $A$ with the corresponding group ring element $\sum_{g\in A}g$, which will again be denoted by $A$. We will be using  the Fourier inversion formula quite frequently.

\begin{lemma}\label{inv} {\em (Inversion Formula)}
Let $G$ be an abelian group of order $N$ and $A=\sum_{g\in G}a_gg\in \C[G]$. Then
 $$a_g=\frac{1}{N}\sum_{\chi\in\widehat G}\chi(A)\chi(g^{-1})$$
for all $g\in G$, where $\widehat{G}$ is the group of complex characters of $G$. Hence if $A, B\in \C[G]$ satisfy $\chi(A)=\chi(B)$ for all $\chi\in \hat G$, then $A=B$.

\end{lemma}

Now set $G=Z_N$, a cyclic group of order $N$ with a generator ${\overline \gamma}$. That is, $Z_N=\{1, {\overline \gamma}, \ldots, {\overline \gamma}^{N-1}\}$. A circulant matrix $M$ in (\ref{cirweigh}) satisfies $MM^{\top}=wI$ if and only if $DD^{(-1)}=w$, where $D$ is the group ring element in $\C[Z_N]$ defined by $D=\sum_{i=0}^{N-1}a_i{\overline \gamma}^i$. Since $a_i=0, \pm 1$, we can write $D=A-B$, where $A=\{{\overline \gamma}^i \mid 0\leq i\leq N-1, a_i=1\}$ and $B=\{{\overline \gamma}^i\mid 0\leq i\leq N-1, a_i=-1\}$. Thus a circulant weighing matrix of order $N$ and weight $w$ is equivalent to a group ring element $A-B$, where $A$ and $B$ are disjoint subsets of $Z_N$, such that
$$(A-B)(A-B)^{(-1)}=w\cdot 1\; {\rm in}\; \C[Z_N].$$

Next we give a short introduction to association schemes. Let $X$ be a finite set. A (symmetric) {\it association scheme} with $d$ classes on $X$ is a partition of
$X\times X$ into subsets $R_0$, $R_1, \ldots , R_d$ (called {\it associate classes} or {\it relations}) such that
\begin{enumerate}
\item $R_0=\{(x,x) \mid  x\in X\}$ (the diagonal relation),
\item $R_i$ is symmetric for $i=1,2,\ldots ,d$,
\item for all $i,j,k$ in $\{0,1,2,\ldots ,d\}$ there is an integer $p_{ij}^k$ such that, for all $(x,y)\in R_k$,
$$|\{z\in X \mid (x,z)\in R_i\; {\rm and}\; (z,y)\in R_j\}|=p_{ij}^k.$$
\end{enumerate}
We denote such an association scheme by $(X, \{R_i\}_{0\leq i\leq d})$. For $i\in \{0,1,\ldots ,d\}$, let $A_i$ be the adjacency matrix of the relation $R_i$, that is, the rows and columns of $A_i$ are both indexed by $X$ and
$$(A_i)_{xy}:=\biggm\{
\begin{array}{c} 1 \quad \mbox{ if }\quad (x,y)\in R_i, \\
                 0 \quad \mbox{ if }\quad (x,y)\notin R_i. \\
\end{array} $$  
The matrices $A_i$ are symmetric $(0,1)$-matrices and
$$A_0=I, \; A_0+A_1+\cdots +A_d=J,$$
where $J$ is the all-1 matrix of size $|X|$ by $|X|$.

By the definition of an association scheme, we have $A_iA_j=\sum_{k=0}^d p_{ij}^kA_k $
for any $i,j\in \{0,1,\ldots ,d\}$; so $A_0,A_1,\ldots , A_d$ form a basis of the commutative algebra generated by $A_0,A_1,\ldots , A_d$ over the reals, which is called the {\it Bose-Mesner algebra} of the association scheme. Moreover this algebra has a unique basis $E_0, E_1, \ldots, E_d$ of primitive idempotents; one of the primitive idempotents is $\frac {1}{|X|}J$.  We may assume that $E_0=\frac {1} {|X|}J$. Let $m_i={\rm rank}\;E_i$. Then
$$m_0=1,\; m_0+m_1+\cdots +m_d=|X|.$$
The numbers $m_0,m_1,\ldots ,m_d$ are called the {\it multiplicities} of the scheme. Since we have two bases of the Bose-Mesner algebra, we consider the transition matrices between them. Define
$P=\left(p_j(i)\right)_{0\le i,j\le d}$ (the {\it first eigenmatrix or character table}) and $Q=\left(q_j(i)\right)_{0\le i,j\le d}$
(the {\it second eigenmatrix}) as the $(d+1)\times (d+1)$ matrices with rows and columns indexed by $0,1,2,\ldots,d$ such that
$$(A_0,A_1, \ldots ,A_d)=(E_0,E_1, \ldots ,E_d)P,$$
and
$$|X|(E_0,E_1, \ldots ,E_d)=(A_0,A_1, \ldots ,A_d)Q.$$
Let $k_j=p_j(0)$, $0\leq j\leq d$. The $k_j$'s  are called {\it valencies} of the scheme.

We call an association scheme $(X,\{R_i\}_{i=0}^d)$, where $X$ is an additively written finite abelian group,  a {\it translation association scheme} or a {\it Schur ring} if  there is a partition $D_0=\{0\},D_1,\ldots,D_d$  of $X$ such that for each $i=0,1,\ldots ,d$,
\begin{equation}\label{eq:ASrelation}
R_i=\{(x,x+y)|\,x\in X, y \in D_i\}.
\end{equation}
Assume that $(X, \{R_i\}_{0\leq i\leq d})$ is a translation association scheme with relations defined in \eqref{eq:ASrelation}.
There is an equivalence relation defined on the character group $\widehat{X}$ of $X$ as follows: $\chi\sim\chi'$ if and only if $\chi(D_i)=\chi'(D_i)$ for all $0\leq i\leq d$. Here $\chi(D)=\sum_{g\in D}\chi(g)$, for any $\chi\in \widehat{X}$ and $D\subseteq X$. Denote by $D_0', D_1',\ldots, D_d'$ the equivalence classes, with $D_0'$ consisting of only the principal character. Define
\[
R_i'=\{(\chi,\chi \chi')|\,\chi\in \widehat{X}, \chi' \in D_i'\}.
\]
Then,
 $(\widehat{X},\{R_i'\}_{i=0}^d)$ also forms a translation association scheme, called the {\it dual} of $(X, \{R_i\}_{0\leq i\leq d})$. The first eigenmatrix of the dual scheme is equal to the second eigenmatrix of the original scheme. We refer the reader to \cite[p.~68]{bcn} for more details.
A translation scheme is called {\it self-dual} if it is isomorphic to its dual.

As an example of translation association schemes, we mention the cyclotomic scheme, which we define below. Let $q$ be a prime power, $N>1$ be a divisor of $q-1$. Let $C_a, 0\leq a\leq N-1$, be the cyclotomic classes of order $N$ of $\F_q$. Assume that $-1\in C_0$. Define $R_0=\{(x,x) \mid  x\in \F_q\}$, and for $a\in \{1,2,\ldots ,N\}$, define $R_a=\{(x,y)\mid x,y\in \F_q, x-y\in C_{a-1}\}$. Then $(\F_q, \{R_a\}_{0\leq a\leq N})$ is an association scheme. This is the so-called {\it cyclotomic association scheme of class $N$ over $\F_q$}. The first eigenmatrix $P$ of the cyclotomic scheme of class $N$ is given by the following $(N+1)$ by $(N+1)$ matrix (with the rows of $P$ arranged in a certain way)
\begin{equation}\label{eigenmatrix}
P=\left(\begin{array}{cccccc}
1& k& k & k &\cdots & k\\
1&\eta_{_{N-1}}  &\eta_0   &\eta_1 & \cdots     &\eta_{_{N-2}}   \\
1&\eta_{_{N-2}}  &\eta_{_{N-1}} & \eta_0& \cdots  &\eta_{_{N-3}} \\
\vdots & & & & \\
1&\eta_0 &\eta_1 &\eta_2& \cdots  &\eta_{_{N-1}}\\
\end{array}\right)
\end{equation}
where $k=\frac{q-1}{N}$ and $\eta_a$, $0\leq a\leq N-1$,  are the Gauss periods of order $N$ defined above. For future use, the submatrix $P_0=(p_j(i))_{1\leq i,j\leq N}$ of $P$ will be called {\it the principal part} of $P$. Note that the cyclotomic scheme $(\F_q, \{R_a\}_{0\leq a\leq N})$ is self-dual.




The rest of the paper is organized as follows. In Section 2, we obtain necessary conditions for Gauss periods to take exactly three rational values. Connections between three-valued Gauss sums and combinatorial structures such as circulant weighing matrices and 3-class association schemes are also developed. In Section 3, we show that in certain cases, the necessary conditions we obtained in Section 2 are also sufficient. Finally in Section 4, we provide five infinite classes of examples where the Gauss periods take exactly three values. Some sporadic examples are also obtained by computer search. From these examples, we obtain circulant weighing matrices and 3-class self-dual association schemes.

\section{Three-valued Gauss periods: Necessary Conditions}

Let $q=p^f$ be a prime power and $N>2$ be a positive integer such that $N\,|\,(q-1)$. Set $k=(q-1)/N$. Let $\F_q$ be the finite field of order $q$, $\gamma$ a fixed primitive element of $\F_q$, and $C_0=\langle \gamma^N\rangle$. Suppose that the Gauss periods $\eta_a=\psi(\gamma^a C_0)$, $a=0,1,\ldots,N-1$, take exactly three distinct rational values $\alpha_1,\alpha_2$, and $\alpha_3$. We will be working with the quotient group $Z_N:=\F_q^\ast/C_0$, a cyclic group of order $N$ with a generator ${\overline \gamma}=\gamma C_0$. For $1\leq i\leq 3$ define subsets $I_i$ of $Z_N$ by $I_i=\{{\overline \gamma}^a\in Z_N\,|\,\eta_a=\alpha_i\}.$

\begin{lemma}\label{lem_gr}  With the above assumptions and notation, we have
\begin{equation}\label{greqn}
(\alpha_1I_1+\alpha_2I_2+\alpha_3I_3)(\alpha_1I_1^{(-1)}+\alpha_2I_2^{(-1)}+\alpha_3I_3^{(-1)})=q\cdot 1-\frac{q-1}{N}Z_{N},
\end{equation}
in the group ring $\Q[Z_N]$.
\end{lemma}

\proof
As above $Z_N$ is the (cyclic) quotient group $\F_q^\ast/C_0$ with a generator ${\overline \gamma}$.  Let $\chi$ be a nontrivial multiplicative character of $\F_q$ whose restriction to $C_0$ is trivial, so that we may view $\chi$ as a character of the quotient group $\F_q^\ast/C_0$. Such a character of $\F_q^\ast/C_0$ will again be denoted by $\chi$ and we have $\chi({\overline \gamma})=\chi(\gamma)$. Note that every nontrivial character of $Z_N:=\F_q^\ast/C_0$ can be obtained in this manner. We have
\begin{eqnarray*}
G_q(\chi)&=&\sum_{i=0}^{N-1}\sum_{x\in C_i}\chi(x)\psi(x)\\
          &=&\eta_0+\chi(\gamma)\eta_1+\cdots +\chi(\gamma^{N-1})\eta_{N-1}\\
          &=&\alpha_1\chi(I_1)+\alpha_2\chi(I_2)+\alpha_3\chi(I_3).
\end{eqnarray*}
Since $\alpha_i$, $1\leq i\leq 3$, are rational integers, we have ${\overline \alpha_i}=\alpha_i$ for $1\leq i\leq 3$. It follows that $$\left(\sum_{i=1}^{3}\alpha_i\chi(I_i)\right)\left(\sum_{i=1}^3\alpha_i\overline {\chi(I_i)}\right)=G_q(\chi)\overline{G_q(\chi)}=q.$$
If $\chi$ is the trivial multiplicative character of $\F_q$, we have $\alpha_1\chi(I_1)+\alpha_2\chi(I_2)+\alpha_3\chi(I_3)=G_q(\chi)=-1$, and $$\left(\sum_{i=1}^{3}\alpha_i\chi(I_i)\right)\left(\sum_{i=1}^3\alpha_i\overline {\chi(I_i)}\right)=G_q(\chi)\overline{G_q(\chi)}=1.$$
The claimed group ring equation now follows from the inversion formula stated in Lemma~\ref{inv}.\qed
\\

Using Lemma~\ref{lem_gr}, we can express the sizes of $I_i$'s in terms of $\alpha_1,\alpha_2,\alpha_3$.

\begin{lemma}\label{rem_sz} Suppose that $\eta_a$, $0\leq a\leq N-1$, take three distinct rational values $\alpha_1, \alpha_2$, and $\alpha_3$. With notation as above and $k=\frac{q-1}{N}$, we have
\begin{align*}
|I_1|=-\frac{\alpha_2 \alpha_3(q-1)+k(q-k+\alpha_2+\alpha_3)}{k(\alpha_1-\alpha_2)(\alpha_3-\alpha_1)},\\
|I_2|=-\frac{\alpha_1\alpha_3(q-1)+k(q-k+\alpha_1+\alpha_3)}{k(\alpha_1-\alpha_2)(\alpha_2-\alpha_3)},\\
|I_3|=-\frac{\alpha_1\alpha_2(q-1)+k(q-k+\alpha_1+\alpha_2)}{k(\alpha_2-\alpha_3)(\alpha_3-\alpha_1)}.
\end{align*}
\end{lemma}

\proof First of all it is clear that $|I_1|+|I_2|+|I_3|=N$. Next we have $\alpha_1|I_1|+\alpha_2|I_2|+\alpha_3|I_3|=\sum_{i=0}^{N-1}\eta_i=-1$. Finally, by comparing the coefficient of the identity element in the two sides of the group ring equation (\ref{greqn}), we get $\alpha_1^2|I_1|+\alpha_2^2|I_2|+\alpha_3^2|I_3|=q-k$. These three equations now uniquely determine $|I_1|,|I_2|$, and $|I_3|$. The proof is complete.
\qed

Next  we derive necessary conditions when the Gauss periods take exactly three values.

\begin{proposition}\label{nece}
Let $q=p^f$ be a prime power and $N>2$ be a positive integer such that $N\,|\,(q-1)$. Assume that the Gauss periods $\eta_a$, $0\leq a\leq N-1$, take exactly three rational values $\alpha_1,\alpha_2,\alpha_3$, say, $\alpha_1-\alpha_2=-tu<0$ and $\alpha_3-\alpha_2=tv>0$ with $t>0$ and $\gcd{(u,v)}=1$. Then $t$ is a power of $p$, and
there exist two positive  integers $r,s$, $0<r,s<N$, such that
\begin{itemize}
\item[(i)] $t(-ur+vs)\equiv -1\,(\mod{N})$;
\item[(ii)] $(N-1)q+t^2(-ur+vs)^2=Nt^2(u^2r+v^2s)$.
\end{itemize}
In particular, $t$ is the largest power of $p$ dividing $G_q(\chi)$ for all nontrivial multiplicative character $\chi$ of $\F_q$ of order dividing $N$.
\end{proposition}

\proof
As before, let $Z_N=\F_q^{\ast}/C_0=\langle {\overline \gamma}\rangle$. So ${\widehat Z_N}=C_0^{\perp}:=\{\chi\mid \chi\in {\widehat \F_q^{\ast}},\;  \chi|_{C_0}=1\}$. We define a function $\v: Z_N\rightarrow \C$ by $\v({\overline \gamma}^a)=\eta_a-\alpha_2.$ In order to simplify notation, we will sometimes write $\v({\overline \gamma}^a)$ simply as $\v(a)$. The Fourier transform of $\v$ is $\widehat{\v}: {\widehat Z_N}\rightarrow \C$, which is defined by
$$\widehat{\v}(\chi)=\frac{1}{\sqrt{N}}\sum_{a=0}^{N-1}\v({\overline \gamma}^a)\chi({\overline \gamma}^a).$$
Computing the Fourier transform of $\v$, we have
\begin{eqnarray*}
\widehat{\v}(\chi)=\left\{
\begin{array}{ll}
\frac{1}{\sqrt{N}}G_q(\chi), &  \mbox{if $\chi$ is a nontrivial character of}\; Z_N;\\
-\frac{1}{\sqrt{N}}-\alpha_2 {\sqrt{N}}, &  \mbox{if $\chi$ is trivial.}
 \end{array}
\right.
\end{eqnarray*}
By assumption $\eta_a$, $0\leq a\leq N-1$, take exactly
three values, we see that $\v(a)\in \{0,-tu,tv\}$.
Note that if $\chi\in {\widehat Z_N}$ is nontrivial, then
\[
G_q(\chi)=\sum_{a=0}^{N-1}\eta_a\chi(\gamma^a)=\sum_{a=0}^{N-1}(\eta_a-\alpha_2)\chi(\gamma^a)=
t(-u\chi(I_1)+v\chi(I_3)),
\]
where
$I_1$ and $I_3$ are defined as before. From the above equation, we see that $t\,|\,G_q(\chi)$ for all nontrivial $\chi\in C_0^{\perp}$. It follows that $t=p^\theta$ for some integer $\theta$.

Let $(C_0^{\perp})^*:=C_0^{\perp}\setminus \{\chi_0\}$ with $\chi_0$ the trivial character. 
Since for any $a, 0\leq a\leq N-1$,
\[
\eta_a=\frac{1}{N}\sum_{\chi\in C_0^\perp}G_q(\chi)\chi^{-1}(\gamma^a),
\]
we have
\begin{equation}\label{delta}
\sigma({\overline \gamma}^a)=\frac{1}{N}\sum_{\chi\in (C_0^\perp)^\ast}G_q(\chi)(\chi^{-1}(\gamma^a)-\chi^{-1}(\gamma^e)),
\end{equation}
where $\chi^{-1}$ is the inverse of $\chi\in C_0^{\perp}$, and $\alpha_2$ is assumed to be equal to $\eta_e$ for some $e$. Let $t'$ be the largest power of $p$ dividing all $G_q(\chi)$, $\chi\in (C_0^\perp)^\ast$. Then (\ref{delta}) implies that $t=t'$ since $\gcd{(N,t')}=1$.

Moreover, by the definition of $\widehat{\v}$ we have
\[
\widehat{\v}(\chi_0)=\frac{1}{\sqrt{N}}\sum_{a=0}^{N-1}\v({\overline \gamma}^a)=
\frac{t(-ur+vs)}{\sqrt{N}},
\]
where $r=|I_1|$, $s=|I_3|$, and $0<r,s<N$. Hence $-\frac{1}{\sqrt{N}}-\frac{N}{\sqrt{N}}\alpha_2=\frac{t(-ur+vs)}{\sqrt{N}}$. It follows that $t(-ur+vs)\equiv -1\,(\mod{N})$.

It is clear from the definition of $\v$ that $\sum_{a=0}^{N-1}\v(a)\overline{\v(a)}=t^2(u^2r+v^2s)$.
On the other hand, we have
\[
\sum_{\chi\in C_0^\perp}\widehat{\v}(\chi)\overline{\widehat{\v}(\chi)}=
\frac{1}{N}\sum_{\chi\in (C_0^\perp)^\ast}G_q(\chi)\overline{G_q(\chi)}+\frac{t^2(-ur+vs)^2}{N}=\frac{1}{N}((N-1)q+t^2(-ur+vs)^2).
\]
It now follows from Parseval's identity that
\[
(N-1)q+t^2(-ur+vs)^2=Nt^2(u^2r+v^2s).
\]
The proof is now complete. \qed

\begin{remark}{\rm
(1) As seen from the proof above, $t$ is the largest power of $p$ dividing all $G_q(\chi)$, $\chi\in (C_0^\perp)^\ast$. By the Stickelberger theorem on the prime ideal factorization of Gauss sums, we have $t=p^{\theta}=p^{d\theta'}$
with $\theta'=\frac{1}{p-1}{\rm min}\{s_p(jk)\mid 1\leq j\leq N-1\}$ and $d=f/f'$, where $f'$ is the order of $p$ modulo $N$ and $s_p(\cdot)$ is the $p$-adic digit sum function.\\

(2) In Section 4, we will show that the two simple necessary conditions in Proposition~\ref{nece} are sometimes also sufficient.}


\end{remark}

\subsection{Circulant Weighing matrices}  Let $q=p^f$ be a prime power, $\gamma$ be a primitive element of $\F_q$,  and $N>1$ be a positive integer such that $N\,|\,\frac{q-1}{p-1}$. In \cite{sw}, it was shown that if the Gauss periods $\eta_a=\psi(\gamma^aC_0)$, $0\leq a\leq N-1$, take exactly two values $\alpha_1$ and $\alpha_2$, then each of the index sets $I_i=\{a\in \Z_N\,|\,\eta_a=\alpha_i\}$, $1\leq i\leq 2$, forms a difference set in $\Z_{N}$, which is a {\it subdifference set} of the Singer difference set. It is natural to ask: if the Gauss periods take exactly three values, what combinatorial structures can we obtain from the index sets $I_1, I_2$ and $I_3$? In this subsection, we will see that under certain conditions, three-valued Gauss periods lead to circulant weighing matrices.

\begin{lemma}\label{AP_size}
Let $q=p^f$ be a prime power and $N>2$ be a positive integer such that $N\,|\,(q-1)$. Assume that the Gauss periods $\eta_a$, $0\leq a\leq N-1$, take exactly three rational values $\alpha_1,\alpha_2,\alpha_3$ which form an arithmetic progression, say, $\alpha_1-\alpha_2=-t<0$ and $\alpha_3-\alpha_2=t>0$. Then
\begin{align*}
&|I_1|=\frac{N(\alpha_2^2+\alpha_2t+k)+2\alpha_2-k+t+1}{2t^2},\;|I_3|=\frac{N(\alpha_2^2-\alpha_2 t+k)+2\alpha_2-k-t+1}{2t^2},\\
&|I_2|=\frac{N(t^2-\alpha_2^2-k)-1-2\alpha_2+k}{t^2},\;\;\; |I_1|-|I_3|=\frac{\alpha_2N+1}{t}.
\end{align*}
Moreover we have
\begin{equation}\label{GReq}
(I_1-I_3)(I_1-I_3)^{(-1)}=\frac{q}{t^2}\cdot 1+\frac{\alpha_2^2N+2\alpha_2-k}{t^2}Z_N
\end{equation}
in $\Q[Z_N]$. In particular, $t$ must be a power of $p$.
\end{lemma}

\proof The fact that $t$ is a power of $p$ follows from Proposition~\ref{nece}. The sizes of $I_1, I_2$ and $I_3$ can be obtained from Lemma~\ref{rem_sz} and the assumptions that $\alpha_1=\alpha_2-t$ and $\alpha_3=\alpha_2+t$. Finally by Lemma~\ref{greqn} and the assumptions that $\alpha_1=\alpha_2-t$ and $\alpha_3=\alpha_2+t$, we have
$$(\alpha_2 Z_N-t(I_1-I_3))(\alpha_2 Z_N-t(I_1-I_3))^{(-1)}=q\cdot 1-\frac{q-1}{N} Z_N,$$
from which~(\ref{GReq}) follows. This completes the proof.
\qed
\vspace{0.1in}

We further consider the question of when $I_1-I_3$ generates a circulant weighing matrix.

\begin{proposition}\label{circu}
Let $q=p^f$ be a prime power and $N>2$ be a positive integer such that $N\,|\,(q-1)$. Assume that the Gauss periods $\eta_a$, $0\leq a\leq N-1$, take exactly three rational values $\alpha_1,\alpha_2,\alpha_3$ which form an arithmetic progression, say, $\alpha_1-\alpha_2=-t<0$ and $\alpha_3-\alpha_2=t>0$. Then
$I_1-I_3$ generates a circulant weighing matrix {\bf CW}$(N, \frac{q}{t^2})$ if and only if $\alpha_2=(\sqrt{q}-1)/N$ and $q$ is a square.
\end{proposition}

\proof Let $q$ be a square and $\alpha_2=(\sqrt{q}-1)/N$. Then $\alpha_2^2N+2\alpha_2-k=0$; in this case (\ref{GReq}) becomes
\[
(I_1-I_3)(I_1-I_3)^{(-1)}=\frac{q}{t^2}\cdot 1,
\]
that is, $I_1-I_3$ generats a circulant weighing matrix of order $q$ and weight $\frac{q}{t^2}$.

Conversely, if $I_1-I_3$ generates a circulant weighing matrix {\bf CW}$(N, \frac{q}{t^2})$, then $\alpha_2^2N+2\alpha_2-k=0$. It follows that $\alpha_2=\frac{\sqrt{q}-1}{N}$ or $\alpha_2=-\frac{1+\sqrt{q}}{N}$. In the latter case, $\sqrt{q}\equiv -1\pmod{N}$, from which we know that
the Gauss periods take only two values~\cite{sw}. Therefore, we must have $\alpha_2=\frac{\sqrt{q}-1}{N}$. Since $\alpha_2$ is rational, we see that $q$ is a square.
\qed


\subsection{Related association schemes} As we remarked in Section 1, when the Gauss periods $\eta_a$, $0\leq a\leq N-1$, take exactly two distinct values, and $-1\in C_0$, then we naturally obtain a strongly regular Cayley graphs defined on $\F_q$ with connection set $C_0$ (which is denoted by ${\rm Cay}(\F_q, C_0)$). Strongly regular graphs are the same objects as 2-class association schemes. We will see in this section that if the Gauss periods take exactly three values, under certain conditions, we obtain 3-class self-dual association schemes. Before stating our main theorem, we give some remarks on translation schemes.

Let $G=\{g_1,\ldots,g_v\}$ be a multiplicative abelian group of order $v$, with character group $\widehat{G}$. Let $\rho:\,G\rightarrow GL_v(\mathbb{C})$ be the regular representation of $G$, namely $\left(\rho(g)\right)_{(h_1,h_2)}=1$ if $h_2=h_1g$, and $=0$ otherwise. Also, for a character $\chi\in \widehat{G}$, let ${\bf v}_\chi:=\frac{1}{\sqrt{v}}(\chi(g_1),\ldots,\chi(g_v))$ and $E_\chi:={\bf v}_\chi^{\top }\cdot \overline{{\bf v}_\chi}$. Then the ${\bf v}_\chi^{\top}$ are the common eigenvectors of $\rho(g)$, $g\in G$, since $\rho(g){\bf v}_\chi^{\top}=\chi(g){\bf v}_\chi^{\top}$. The $E_\chi$'s are the primitive idempotents of the algebra $\mathcal{A}:=\langle \rho(g):g\in G\rangle\cong \mathbb{C}[G]$, as can be easily checked by using the orthogonal relations of characters. Moreover, by using the fact that $\left(E_{\chi}\right)_{(g,h)}=\frac{1}{v}\chi(gh^{-1})$, we have
\begin{align}\label{eqn_prod}
(vE_\chi)\circ (vE_{\chi'})=(vE_{\chi\chi'}).
\end{align}

Now assume that $D_0,D_1,\ldots ,D_d$ form a partition of $G$ which yields a translation scheme, and its dual scheme is given by the following partition of $\widehat{G}$: $D_0',D_1',\ldots,D_d'$. Write $A_i=\rho(D_i)$, $0\le i\le d$, and let $E_0,E_1,\ldots,E_d$ be the primitive idempotents of the Bose-Mesner algebra $\mathcal{A}:=\langle A_0,A_1,\ldots,A_d\rangle$ with respect to the matrix multiplication. We have $E_i=\sum_{\chi\in D_i'}E_\chi$ with respect to a proper ordering of the $E_i$'s. 

Similarly, if we use $\rho'$ for the regular representation of $\widehat{G}$, and write $A_i'=\rho'(D_i')$, then $A_0',\ldots,A_d'$ span $\widehat{\mathcal{A}}$,  the  Bose-Mesner algebra of the dual scheme.

Let $\Psi$ be the linear map from $\widehat{\mathcal{A}}$ to $\mathcal{A}$ that maps $A_i'$ to $vE_i$, $0\le i\le d$. It follows from \eqref{eqn_prod} that $\Psi$ is an algebra isomorphism from $(\widehat{\mathcal{A}},+,\cdot)$ to $(\mathcal{A},+,\circ)$.  An easy corollary is that, $\Psi$ maps the idempotents of $(\widehat{\mathcal{A}},+,\cdot)$ to those of $(\mathcal{A},+,\circ)$, namely, the $A_i$'s. 

\begin{theorem}\label{selfasso}
Let $q=p^f$ be a prime power and $N>2$ be a positive integer such that $N\,|\,(q-1)$. Assume that $-1\in C_0$ and the Gauss periods $\eta_a$, $0\leq a\leq N-1$, take exactly three rational values $\alpha_1,\alpha_2,\alpha_3$, say, $\alpha_1-\alpha_2=-tu<0$ and $\alpha_3-\alpha_2=tv>0$ with $t>0$. Let
\[
R_0=\{0\},\, R_1=\bigcup_{i\in I_1}C_i, \,R_2=\bigcup_{i\in I_2}C_i, \, R_3=\bigcup_{i\in I_3}C_i,
\]
where $I_i=\{a\in \Z_N\,|\,\eta_a=\alpha_i\}$ for $i=1,2,3$. If $|I_1|=1$ or $|I_3|=1$, then $(\F_q,\{{\mathcal R}_i\}_{i=0}^{3})$ is a self-dual three-class association scheme. (Here for $0\leq i\leq 3$, $(x,y)\in {\mathcal R}_i$ if and only if $x-y\in R_i$.)
\end{theorem}
\proof 
Let $A_0,A_1,\ldots,A_N$ and $E_0,E_1,\ldots,E_N$ be the first and the
second standard bases of the Bose-Mesner algebra ${\mathcal A}$ of
the cyclotomic scheme of class $N$ of $\F_q$. We may assume that
 the cyclic permutation $\sigma=(1,2,\ldots,N)$ is an algebraic automorphism
of the association scheme, namely, the linear map  that maps $A_i\mapsto A_{\sigma(i)}$, $0\le i\le d$, is an automorphism of the Bose-Mesner algebra with respect to both the matrix multiplication and the Schur product.  Notice that ${\mathcal A}$ consists of symmetric matrices.

In what follows  we use the notation $E_S=\sum_{s\in S}E_s$, $A_S=\sum_{s\in S} A_s$ for any $S\subseteq \{1,2,\ldots,N\}$. Let $P=(p_j(i))$ and $Q=(q_j(i))$ be the first and second eigenmatrix of the cyclotomic scheme, respectively. With a proper ordering of the  $E_i$'s, we have $P=Q$, and the principal part of $P$ is symmetric. The principal part of $P$ has only three distinct rational entries, namely, $\alpha_1,\alpha_2,\alpha_3$.

Since the Gauss periods have three values from $\{\alpha_1,\alpha_2,\alpha_3\}$, we have \[A_1=kE_0+\alpha_1 E_{L_1}+\alpha_2 E_{L_2}+\alpha_3 E_{L_3},\]
where the $L_i$'s form a partition of $\{1,2,\ldots,N\}$ (and they come from the $I_i$'s in the statement of the theorem). Since the cyclotomic scheme is self-dual, by the algebra isomorphism $\Psi$ described right before the statement of this theorem,
we have 
\[E_1=q^{-1}(kA_0+\alpha_1 A_{M_1}+\alpha_2 A_{M_2}+\alpha_3 A_{M_3}),\] 
where  the $M_i$'s form a partition of $\{1,2,\ldots,N\}$, and $|M_i|=|L_i|$, for all $i=1,2,3$.

Assume now that $|L_1|=1$, i.e., $L_1=\{\ell\}$ for some $\ell\in \{1,2,\ldots,N\}$. We have $\alpha_1=p_1(\ell)$. Consider the vector space ${\mathcal B}$ spanned by
$A_0,A_1,E_0,E_\ell$. Noting that $\alpha_2\not=\alpha_3$, it follows from
\begin{align*}
A_0&\,=E_0+E_\ell+E_{L_2}+E_{L_3},\\
A_1&\,=kE_0+\alpha_1E_\ell+\alpha_2 E_{L_2}+\alpha_3 E_{L_3},
\end{align*}
that ${\mathcal B}=\langle E_0,E_\ell,E_{L_2},E_{L_3}\rangle$. In particular,
${\mathcal B}$ is closed with respect to the matrix multiplication. 

Since $\sigma$
is an algebraic automorphism of ${\mathcal A}$, we have
$E_\ell=q^{-1}(kA_0+\alpha_1 A_{M_1'}+\alpha_2 A_{M_2'}+\alpha_3 A_{M_3'})$,  where $M_i'=\sigma^{\ell-1}(M_i)$. It follows from
$|M_i'|=|M_i|=|L_i|$ that $M_1'=\{m\}$ for some $m\in \{1,2,\ldots,N\}$.
So, $E_\ell=q^{-1}(kA_0+\alpha_1 A_{m}+\alpha_2 A_{M_2'}+\alpha_3 A_{M_3'})$.
On the other hand, we have  
\[
qE_\ell\circ A_1= q_{\ell} (1) A_1=p_{\ell}(1)A_1=p_{1}(\ell)A_1=\alpha_1A_1,
\]
It follows that $m=1$. Together with $E_0=q^{-1}(A_0+A_1+A_{M_2'}+A_{M_3'})$,
we see that ${\mathcal B}=\langle A_0,A_1,A_{M_2'},A_{M_3'}\rangle$. In particular,
${\mathcal B}$ is closed with respect to the Schur product.

Since $A_0,E_0\in {\mathcal B}$ and ${\mathcal B}$ is symmetric, we conclude that
$(\F_q,\{{\mathcal R}_i\}_{i=0}^{3})$ is a self-dual association scheme with
${\mathcal B}$ as its Bose-Mesner algebra.
\qed 

\begin{remark}{\rm
We comment that the condition $|I_1|=1$ (or $|I_3|=1$) in Theorem~\ref{selfasso} is needed. Below is an example in which the Gauss periods take three values, but the partition of $\Z_N$ by $I_1, I_2$ and $I_3$ does not yield a three-class association scheme. Let $q=11^3$, $N=19$, $I_1=\{ 0, 2, 3, 4, 5, 6, 9, 14, 16, 17\}$, $I_2=\{ 8, 10, 12, 13, 15, 18\}$,  and $I_3=\{ 1, 7, 11\}$.
In this case,
$\psi(\gamma^a C_0)$, $a=0,1,\ldots,N-1$, take the values $-7,4$, and $15$ according as $a\in I_i$, $1\leq i\leq 3$,
but the partition $I_1,I_2,I_3$ of $\Z_N$ does not yield a three-class association scheme. }
\end{remark}

\section{Sufficient conditions for Gauss periods to take exactly three values}

In this section, we consider the question  when the necessary conditions obtained in Proposition~\ref{nece} are also sufficient. We pay special attention to the case where either $u=1$ or $v=1$. Here we are using the notation of Proposition~\ref{nece}. (Many  examples given in Section~\ref{ex3-v} fall into this case.)  Furthermore, we show that the partition of $\Z_N$ by $I_1,I_2$, and $I_3$ yields a $3$-class association scheme if  $u=|I_3|=1$ or $v=|I_1|=1$.

\subsection{Sufficient conditions for Gauss periods to take three values}
In this subsection, we give sufficient conditions for Gauss periods to take exactly three distinct values.
First, we give a general sufficient condition. Below we use $\N$ to denote the set of positive integers.

\begin{proposition}\label{suff0}
Let $q=p^f$ be a prime power, $N>2$ be an integer such that $N\,|\,(q-1)$, and $C_0=\langle \gamma^N\rangle$, where $\gamma$ is a fixed primitive element of $\F_q$. Assume that there are four positive integers $u,v,r,s$ such that
\begin{itemize}
\item[(i)] $t(-ur+vs)\equiv -1\,(\mod{N})$;
\item[(ii)] $(N-1)q+t^2(-ur+vs)^2=Nt^2(u^2r+v^2s)$,
\end{itemize}
where $t$ is the largest power of $p$ dividing all $G_q(\chi)$, $\chi\in (C_0^\perp)^\ast=(C_0)^{\perp}\setminus\{\chi_0\}$.
If all nonnegative solutions $(t_x)_{x\in \Z\setminus\{0\}}$ to the following system of equations

$$\begin{cases}\label{eq21}
\sum_{x\in \N}x(x-1)t_x+\sum_{x\in \N}x(x+1)t_{-x}=u(u+1)r+v(v-1)s\\
\sum_{x\in \N}x(x+1)t_x+\sum_{x\in \N}x(x-1)t_{-x}=u(u-1)r+v(v+1)s\\
\end{cases}$$

\noindent satisfy $t_{x}\not=0$ if $x=i_1$ or $i_2$, $t_x=0$ for all $x\neq i_1, i_2$, and $t_{i_1}+t_{i_2}<N$, where $i_1, i_2$ are two distinct integers,
then the Gauss periods $\eta_a=\psi(\gamma^a C_0)$, $0\leq a\leq N-1$, take exactly
three distinct values.
\end{proposition}
\proof
Let $y=\frac{-t (-ur+vs)-1}{N}$. Define a map $\tau:\Z_N\to \C$ by
\[
\tau(a)=\frac{\psi(\gamma^a C_0)-y}{t}.
\]
Since
\[
\psi(\gamma^a C_0)+\frac{1}{N}=\frac{1}{N}\sum_{\chi\in (C_0^\perp)^\ast}
G_q(\chi)\chi^{-1}(\gamma^a),
\]
by $t\,|\,G_q(\chi)$ and assumption (i), we see that
$\tau(a)\in \Z$, that is,  $\tau$ is integer-valued.
Computing the Fourier transform of $\tau$, we have
\begin{eqnarray*}
\widehat{\tau}(\chi)=\frac{1}{\sqrt{N}}\sum_{a\in \Z_N}\tau(a)\chi(\gamma^a)=\left\{
\begin{array}{ll}
\frac{1}{t \sqrt{N}}G_q(\chi),&  \mbox{if $\chi$ is nontrivial;}\\
\frac{-ur+vs}{\sqrt{N}},&  \mbox{if $\chi$ is trivial.}
 \end{array}
\right.
\end{eqnarray*}
It follows from Parseval's identity that
\[
\sum_{a\in \Z_N}\tau(a)^2=
\sum_{\chi\in C_0^\perp}
\widehat{\tau}(\chi)\overline{\widehat{\tau}(\chi)}=(N-1)\frac{q}{Nt^2}+\frac{(-ur+vs)^2}{N}.
\]
By assumption (ii), we have

\begin{equation}\label{2ndmoment}
\sum_{a\in \Z_N}\tau(a)^2=u^2r+v^2s.
\end{equation}
On the other hand, we have

\begin{equation}\label{1stmoment}
\sum_{a\in \Z_N}\tau(a)=-ur+vs.
\end{equation}
Eqs.~(\ref{2ndmoment}) and (\ref{1stmoment}) can be rewritten as
\[
\sum_{x\in \N}x^2t_x+\sum_{x\in \N}x^2t_{-x}=u^2r+v^2s\, \, \mbox{ and }\, \,
\sum_{x\in \N}xt_x-\sum_{x\in \N}xt_{-x}=-ur+vs,
\]
where $t_x=|\{a\in \Z_N\,|\,\tau(a)=x\}|$, $x\in \N$. It follows that
\begin{equation}\label{subtracting}
\sum_{x\in \N}x(x-1)t_x+\sum_{x\in \N}x(x+1)t_{-x}=u(u+1)r+v(v-1)s
\end{equation}
and
\begin{equation}\label{adding}
\sum_{x\in \N}x(x+1)t_x+\sum_{x\in \N}x(x-1)t_{-x}=u(u-1)r+v(v+1)s.
\end{equation}
By assumption, the nonnegative solutions $(t_x)_{x\in \Z\setminus\{0\}}$ to the above system of equations all satisfy
$t_{x}\not=0$ when $x=i_1$ or $i_2$ and $t_x=0$ for all $x\neq i_1, i_2$.  This implies that $\tau(a)\in \{0,i_1,i_2\}$ for all $a\in \Z_N$. Consequently $\eta_a=\psi(\gamma^a C_0)$, $0\leq a\leq N-1$, take exactly three distinct values since $t_{i_1}+t_{i_2}<N$. The proof is complete.
\qed
\vspace{0.3cm}

As an immediate corollary, we have the following.
\begin{corollary}\label{co:self}
Let $q=p^f$ be a prime power, $N>2$ be an integer such that $N\,|\,(q-1)$, and $C_0=\langle \gamma^N\rangle$, where $\gamma$ is a fixed primitive element of $\F_q$. Assume that there are four positive integers $u,v,r,s$ satisfying
\begin{itemize}
\item[(i)] $t(-ur+vs)\equiv -1\,(\mod{N})$;
\item[(ii)] $(N-1)q+t^2(-ur+vs)^2=Nt^2(u^2r+v^2s)$,
\end{itemize}
where $t$ is the largest power of $p$ dividing all $G_q(\chi)$, $\chi\in (C_0^\perp)^\ast=(C_0)^{\perp}\setminus\{\chi_0\}$. If $u=v=r=1$ and $s+1<N$, or $u=v=s=1$ and $r+1<N$, then $\eta_a=\psi(\gamma^a C_0)$, $0\leq a\leq N-1$, take exactly
three distinct values; in this case, the three values taken by $\eta_a$ form an arithmetic progression.
\end{corollary}
\proof
We assume that  $u=v=s=1$. (The case where $u=v=r=1$ is similar.)
In this case, (\ref{adding}) is reduced to
\[
\sum_{x\in \N}x(x+1)t_x+\sum_{x\in \N\setminus\{1\}}x(x-1)t_{-x}=2.
\]
The nonnegative solutions $(t_x)_{x\in \Z\setminus\{0\}}$ to the system of equations

$$\begin{cases}
\sum_{x\in \N\setminus\{1\}}x(x-1)t_x+\sum_{x\in \N}x(x+1)t_{-x}=2r\\
\sum_{x\in \N}x(x+1)t_x+\sum_{x\in \N\setminus\{1\}}x(x-1)t_{-x}=2.\\
\end{cases}$$

\noindent must satisfy $t_{1}=1$, $t_{-1}=r$ and $t_x=0$ for all other $x$, or
$t_{-2}=1$, $t_{-1}=r-3$ and $t_x=0$ for all other $x$.  It follows that $\tau(a)\in\{0,-1,1\}$  or $\tau(a)\in\{0,-1,-2\}$ for all $a\in \Z_N$. Consequently $\eta_a$, $a\in \Z_N$, take exactly three distinct values since $s+1<N$. The proof of the corollary is complete.
\qed
\vspace{0.3cm}

The conditions $u=v=r=1$ and $s+1<N$ in the above corollary are quite restrictive. Below we consider more general situations where we can still guarantee that the Gauss periods take only three values. We start with the following lemma.

\begin{lemma}\label{le:size}
Let $q=p^f$ be a prime power, $N>1$ be an integer such that $N\,|\,(q-1)$, and $C_0=\langle \gamma^N\rangle$, where $\gamma$ is a fixed primitive element of $\F_q$. Assume that $\eta_a=\psi(\gamma^a C_0)$, $0\leq a\leq N-1$, take exactly $\ell$ distinct values, say, $\alpha_1,\alpha_2,\ldots,\alpha_\ell$. Let $I_i=\{a\in \Z_N\,|\,\eta_a=\alpha_i\}$ for $1\leq i\leq \ell$. Then each $I_i$ is invariant under the multiplication by $p$. Moreover, assume that $m:=\gcd\{\ord_{n}(p)\,|\,n>1\, \mbox{and $n$ divides $N$}\}\ge 2$. Then there exists a unique $i_0$, $1\leq i_0\leq \ell$, such that $|I_{i_0}|\equiv 1\,(\mod{m})$ and $|I_i|\equiv 0\,(\mod{m})$ for all $i\neq i_0$.
\end{lemma}

\proof Since $\Tr_{q/p}(x^p)=\Tr_{q/p}(x)$ for $x\in \F_q$, we have $\eta_{pa}=\eta_a$ for all $a\in \Z_N$. It follows that each $I_i$ is invariant under the multiplication by $p$. Note that under the multiplication by $p$ (i.e., under the map $x\mapsto px$, $x\in \Z_N$), $0$ forms a singleton orbit, all other orbits have sizes divisible by $m$. The second conclusion of the lemma follows. This completes the proof of the lemma.
\qed

\begin{theorem}\label{suff}
Let $q=p^f$ be a prime power, $N>2$ be an integer such that $N\,|\,(q-1)$, and $C_0=\langle \gamma^N\rangle$, where $\gamma$ is a fixed primitive element of $\F_q$. Assume that there are four positive integers $u,v,r,s$ such that
\begin{itemize}
\item[(i)] $t(-ur+vs)\equiv -1\,(\mod{N})$;
\item[(ii)] $(N-1)q+t^2(-ur+vs)^2=Nt^2(u^2r+v^2s)$,
\end{itemize}
where $t$ is the largest power of $p$ dividing all $G_q(\chi)$, $\chi\in (C_0^\perp)^\ast=(C_0)^{\perp}\setminus\{\chi_0\}$. Let $m=\gcd\{\ord_{n}(p)\,|\,n>1,\mbox{and $n$ divides $N$}\}$ and  assume that $m\ge 2$. If one of the following conditions holds,
\begin{itemize}
\item[(1)] $u=s=1$, $v(v+1)< 2m$, and $r+1<N$;
\item[(2)] $u=s=1$, $v(v+1)=2m$, and $r+v^2<N$;
\item[(3)] $v=r=1$, $u(u+1)< 2m$, and $s+1<N$;
\item[(4)] $v=r=1$, $u(u+1)= 2m$, and $s+u^2<N$;
\item[(5)] $u=v=1$, $s=m$, and $r+m<N$;
\item[(6)] $u=v=1$, $r=m$, and
$s+m<N$,
\end{itemize}
then $\eta_a=\psi(\gamma^a C_0)$, $0\leq a\leq N-1$, take exactly three values.
\end{theorem}

\proof
First we note that by Lemma~\ref{le:size}, the $t_x$, $x\in \Z\setminus\{0\}$, in Eqs.~(\ref{subtracting}) and (\ref{adding}) satisfy that $t_x\equiv 1\,(\mod{m})$ for at most one $x$ and
$m\,|\,t_x$ for all other $x$.

We consider Cases  (1) and (2) where $u=s=1$. (For Cases (3) and (4), the claims can be proved in a similar way. We omit the proof.) In these cases, (\ref{adding}) is reduced to
\begin{equation}\label{eq3}
\sum_{x\in \N}x(x+1)t_x+\sum_{x\in \N\setminus\{1\}}x(x-1)t_{-x}=v(v+1).
\end{equation}

{\bf (1).} If $v(v+1)< 2m$, noting the divisibility conditions on the $t_x$'s, we see that the nonnegative solutions $(t_x)_{x\in \Z\setminus\{0\}}$ to the following system

$$\begin{cases}
\sum_{x\in \N\setminus\{1\}}x(x-1)t_x+\sum_{x\in \N}x(x+1)t_{-x}=2r+v(v-1)\\
\sum_{x\in \N}x(x+1)t_x+\sum_{x\in \N\setminus\{1\}}x(x-1)t_{-x}=v(v+1).\\
\end{cases}$$

\noindent must satisfy $t_v=1$, $t_{-1}=r$ and $t_x=0$ for all other $x$, or
$t_{-(v+1)}=1$, $t_{-1}=r-2v-1$, and $t_x=0$ for all other $x$. It follows that $\tau(a)\in\{0,-1,v\}$ or $\tau(a)\in\{0,-1,-v-1\}$ for all $a\in \Z_N$. Therefore $\eta_a$, $a\in \Z_N$, take exactly three values since $r+1<N$.

{\bf (2).} If $v(v+1)=2m$, the above system has further nonnegative solutions
$t_1=m$, $t_{-1}=r+m-v$ and $t_x=0$ for all other $x$, or $t_{-2}=m$,  $t_{-1}=r-2m-v$, and $t_x=0$ for all other $x$.  So $\tau(a)\in\{0,-1,1\}$ or  $\tau(a)\in\{0,-1,-2\}$ for all $a\in \Z_N$. It follows that $\eta_a$, $a\in \Z_N$, take exactly three values since $r+v^2<N$.

Next, we consider the case where $u=v=1$ and $s=m$ or $r=m$.

{\bf (5).} We assume that $s=m$. (The case where $r=m$ can be handled similarly). In this case, (\ref{adding}) is reduced to
\[
\sum_{x\in \N}x(x+1)t_x+\sum_{x\in \N\setminus\{1\}}x(x-1)t_{-x}=2m.
\]

If $2m\neq \ell(\ell+1)$ for all $\ell\in \Z$, then the nonegative solutions $(t_x)_{x\in \Z\setminus\{0\}}$ to the following system

$$\begin{cases}
\sum_{x\in \N\setminus\{1\}}x(x-1)t_x+\sum_{x\in \N}x(x+1)t_{-x}=2r\\
\sum_{x\in \N}x(x+1)t_x+\sum_{x\in \N\setminus\{1\}}x(x-1)t_{-x}=2m.\\
\end{cases}$$

\noindent must satisfy $t_1=m$,  $t_{-1}=r$, and $t_x=0$ for all other $x$, or $t_{-2}=m$, $t_{-1}=r-3m$, and $t_x=0$ for all other $x$.  It follows that $\tau(a)\in\{0,1,-1\}$ or  $\tau(a)\in\{0,-1,-2\}$ for $a\in \Z_N$. Therefore $\eta_a$, $a\in \Z_N$, take exactly three values.

If $2m$ can be written as $2m=\ell(\ell+1)$ for some positive integer $\ell$, then the above system has further nonnegative solutions $t_\ell=1$, $t_{-1}=r-\ell(\ell-1)/2$, and $t_x=0$ for other $x$, or $t_{-\ell-1}=1$,  $t_{-1}=r-(\ell+1)(\ell+2)/2$, and $t_x=0$ for other $x$. Again we have $\tau(a)\in\{0,\ell,-1\}$ or $\tau(a)\in\{0,-1,-\ell-1\}$ for $a\in \Z_N$.
\qed

\section{Examples of three-valued Gauss periods and related weighing matrices and association schemes}\label{ex3-v}

In this section, we give examples of three-valued Gauss periods. These examples often lead to interesting combinatorial structures such as circulant weighing matrices and association schemes.

As a preparation, we consider a group ring version of the Hasse-Davenport Theorem.

\begin{theorem}\label{thm:lift}{\em (\cite[Theorem~11.5.2]{BEW97})}
Let $\chi$ be a nonprincipal multiplicative character of $\F_q=\F_{p^f}$ and
let $\chi'$ be the lifted character of $\chi$ to the extension field $\F_{q'}=\F_{p^{fe}}$, that is, $\chi'(\alpha):=\chi(\Norm_{q'/q}(\alpha))$ for any $\alpha\in \F_{q'}^*$.
Then, it holds that
\[
G_{q'}(\chi')=(-1)^{e-1}(G_{q}(\chi))^e.
\]
\end{theorem}
Let $\chi$ be a multiplicative character of $\F_q$ of order $N>1$, $\gamma$ a primitive element of $\F_q$, and $C_0=\langle \gamma^N\rangle$. As we saw in the proof of Lemma~\ref{greqn}, we have
$$G_q(\chi)=\eta_0 + \eta_1 \chi(\gamma) + \cdots + \eta_{N-1}\chi(\gamma)^{N-1},$$
where $\eta_a=\psi(C_a^{(N,q)})$ for $0\leq a\leq N-1$. This motivated us to define the following group ring element
\[
g_{F,N}=\sum_{a\in\Z_N}\eta_a[a]\in\C[\Z_N],
\]
where $F=\F_q$. (See \cite{grtf}.) Let $E$ be the finite field with $q^e$ elements, $e>1$ a positive integer. Then it follows from Theorem~\ref{thm:lift} that
\begin{equation}\label{grHD}
g_{E,N}=(-1)^{e-1}g_{F,N}^e.
\end{equation}
The advantage of this group ring version of the Hasse-Davenport theorem is that starting with a pair of small $(q,N)$ with $N|(q-1)$ we are able to determine the Gauss periods corresponding to the subgroup of index $N$ of $\F_{q^e}^{\ast}$ efficiently.

\subsection{Examples from a conic}\label{quadric}
Let $p$ be a prime, $f$ a positive integer,  $F=\F_{p^{3f}}$, and $E=\F_{p^{3fe}}$ with $e>1$. Let $\gamma$ and $\omega$ be primitive elements of $F$ and $E$ respectively such that $\gamma=\Norm_{E/F}(\omega)$. Let $N=\frac{p^{3f}-1}{p^f-1}$. Then $C_0^{(N,F)}=\F_{p^f}^\ast < F^*=\F_{p^{3f}}^*$, and the Gauss periods $\eta_a=\psi(\gamma^aC_0^{(N,F)})=p^f-1$ if $\Tr_{F/L}(\gamma^a)=0$ and $-1$ otherwise, where $L=\F_{p^f}$. Denote by
\[
S:=\{i\in\Z_N:\,\Tr_{F/L}(\gamma^i)=0\}.
\]
Then $|S|= p^f+1$, and $g_{F,N}=p^f S - \Z_N$. As in \cite{fka}, we identify the points of the projective plane $PG(2, p^f)$ with the elements of $\Z_N$. Then $S$ represents a line of $PG(2, p^f)$, and is the well-known Singer difference set in $\Z_N$; see \cite{pott} for instance.

Now set $e=2$. Then by (\ref{grHD}), we have $$g_{E,N}=-(p^fS-\Z_N)^2=-p^{2f}S^2+(p^{2f}+p^f-1)\Z_N.$$ Note that here $g_{E,N}=\sum_{a\in \Z_N}\psi'(\omega^aC_0^{(N,E)}) [a] \in \C[\Z_N]$, $\psi'$ is the canonical additive character of $E$. In order to know how many values the Gauss periods $\psi'(\omega^aC_0^{(N,E)})$, $0\leq a\leq N-1$, take, it suffices to compute $S^2$ in the group ring $\C[\Z_N]$. For any $a\in\Z_N$, the coefficient of $[a]$ in $S^2$ is equal to the size of
\[
\{i\in\Z_N:\,\Tr_{F/L}(\ga^{-i})=0,\Tr_{F/L}(\ga^{i+a})=0\}={\mathcal Q}\cap (S-a),
\]
where ${\mathcal Q}=\{i\in\Z_N:\,\Tr_{F/L}(\ga^{-i})=0\}$ and $S-a=\{x-a\mid x\in S\}$. Since ${\mathcal Q}$ is a conic in $PG(2,p^f)$ (cf. \cite{jv}) and $S-a$ is a line of $PG(2,p^f)$, we have $|{\mathcal Q}\cap (S-a)|=0,1$ or 2, according as $S-a$ is passant, tangent or secant. It follows that the Gauss periods $\psi(\omega^aC_0^{(N,E)})$, $0\leq a\leq N-1$, take three values $\alpha_1=p^{2f}+p^f-1$, $\alpha_2=p^f-1$, and $\alpha_3=-p^{2f}+p^f-1$, which form an arithmetic progression with common difference $t=p^{2f}$. Here $|E|=q^{6f}$ and $\alpha_2=p^f-1= \frac{\sqrt{p^{6f}} -1}{N}$. So by Proposition~\ref{circu} we obtain a {\bf CW}$(p^{2f+p^f+1}, p^{2f})$. We remark that the circulant weighing matrix {\bf CW}$(p^{2f+p^f+1}, p^{2f})$ obtained here is not new (cf. \cite{ww}), but the connection with three-vauled Gauss periods is new.

Note that with the same notation as above, in the special case where $p=2$,  the authors of \cite{fka} already showed that the Cayley graph $Cay(\F_{q},C_0^{(N,q)})$, with $q=2^{6f}$ and $N=(2^{3f}-1)/(2^f-1)$, has three restricted eigenvalues $-2^{2f}+2^f-1,2^f-1,2^{2f}+2^f-1$, and $\{(x,y)\in \F_q\times \F_q\mid x-y\in C_0^{(N,q)}\}$ is a relation in a three-class association scheme, see \cite[p.~1210]{fka}.




\subsection{More examples from two-valued Gauss periods}\label{twoGauss} Let $p$ be a prime, $f\geq 1$ and $e>1$ be integers, and $F=\F_{p^f}$, $E=\F_{p^{fe}}$. Assume that $k|(p^f-1)$. Then certainly $k|(p^{fe}-1)$. Let $N=(p^f-1)/k$ and $N'=(p^{fe}-1)/k$. Then $C_0^{(N,F)}=C_0^{(N',E)}.$ This can be seen as follows. Let $\omega$ and $\gamma$ be primitive elements of $E$ and $F$, respectively, such that $\gamma=\omega^{\frac{p^{fe}-1}{p^f-1}}$. Then $C_0^{(N,F)}=\langle \gamma^N\rangle=\langle \omega^{\frac{(p^{fe}-1)N}{p^f-1}}\rangle=\langle \omega^{N'} \rangle=C_0^{(N', E)}.$

Assume that the Gauss periods $\eta_a=\psi(\gamma^aC_0^{(N,F)}), 0\leq a\leq N-1$, take exactly two distinct values $\alpha_1$ and $\alpha_2$ according as $a\in S$ or not for some $S\subseteq \Z_{N}$. Let $\psi'$ be the canonical additive character of $E$. Then, we have
\begin{align*}
\psi'(\omega^aC_0^{(N', E)})&=\sum_{x\in C_0^{(N,p^{f})}}\xi_p^{{\tr_{p^{f}/p}(x\cdot (\tr_{E/F}(\omega^a)))}}
=\psi(\tr_{E/F}(\omega^a) C_0^{(N,p^{f})})\\
&=\left\{
\begin{array}{ll}
k,&  \mbox{if $\tr_{E/F}(\omega^a)=0$;}\\
\alpha_1,&  \mbox{if $\tr_{E/F}(\omega^a)=\gamma^b$ and $b\in S$;}\\
\alpha_2,&  \mbox{if $\tr_{E/F}(\omega^a)=\gamma^b$ and $b\in \Z_{N}\setminus S$. }
 \end{array}
\right.
\end{align*}
That is, the Gauss periods $\psi'(\omega^aC_0^{(N', E)})$, $0\leq a\leq N'-1$, take three distinct values $k, \alpha_1$ and $\alpha_2$. Furthermore, it is routine to check that
$C_0^{(N,F)}$, $\F^\ast\setminus C_0^{(N,F)}$, $E^\ast\setminus F^\ast$ give a three-class association  scheme.


\subsection{Examples from union of 1-dimensional subspaces}\label{primaryEX}

Let $q\equiv 1\pmod 3$ and $\gamma$ an element of order $k=3(q-1)$ in $\F_{q^3}$, and set $N=\frac{q^3-1}{k}$. Then the degree of the minimal polynomial of $\gamma$ over $\F_q$ is equal to $ord_{k}(q)$. Assume that $ord_{k}(q)=3$. Then $1,\gamma,\gamma^2$ are linearly independent over $\F_q$, and it follows that $C_0^{(N,q^3)}=\langle\gamma^N\rangle=\{\lambda\cdot 1 \mid  \lambda\in \F_q^*\}\cup \{\lambda\cdot \gamma\mid  \lambda\in \F_q^*\}\cup \{\lambda\cdot \gamma^2 \mid  \lambda\in \F_q^*\}$. For any nontrivial additive character $\psi'$ of $\F_{q^3}$ we have

\begin{align*}
\psi'(C_0^{(N,q^3)})=\left\{
\begin{array}{ll}
-3, & \mbox{if $\psi'|_{\F_q}$, $\psi'|_{\F_q\gamma}$, $\psi'|_{\F_q\gamma^2}$ are all nontrivial;}\\
-3+q, &\mbox{if exactly one of $\psi'|_{\F_q}$, $\psi'|_{\F_q\gamma}$, $\psi'|_{\F_q\gamma^2}$ is trivial;}\\
-3+2q, &\mbox{if exactly two of $\psi'|_{\F_q}$, $\psi'|_{\F_q\gamma}$, $\psi'|_{\F_q\gamma^2}$ are trivial.}
\end{array}
\right.
\end{align*}
Therefore the Gauss periods $\eta_a$, $0\leq a\leq N-1$, of $\F_{q^3}$ take three values $\alpha_1=-3, \alpha_2=-3+q, \alpha_3=-3+2q$, which form an arithmetic progression with common difference $t=q$. By Lemma~\ref{AP_size}, we have
$$|I_1|=\frac{(q-1)^2}{3},   |I_2|=q-1,  |I_3|=1.$$
Since $|I_3|=1$, by Theorem~\ref{selfasso}, the subsets $\cup_{i\in I_j}C_i^{(N,q^3)}$, $j=1,2,3$, give a 3-class self-dual association scheme. Note that with assumptons as above, $\alpha_2^2 N+2\alpha_2-k=0$ if and only if $(q,N)=(4,7)$. Therefore we obtain a {\bf CW}$(7,4)$ in the case when $(q,N)=(4,7)$, and we do not obtain circulant weighing matrices in other cases.

\subsection{Examples from products of subfields}\label{sec:subf}
Let $e,f$ be two positive integers such that $e/\gcd{(e,f)}=3$ and  let $q=p^{\lcm{(e,f)}}=p^{3f}$. Let $C_0^{(N,q)}$ be the subgroup of $\F_{q}^\ast$ generated by $\F_{p^{e}}^\ast$ and $\F_{p^{f}}^\ast$. Then
$$|C_0^{(N,q)}|=(p^{e}-1)(p^{f}-1)/(p^\ell-1),$$
where $\ell=\gcd{(e,f)}$ and $N=\frac{(p^{3f}-1)(p^\ell-1)}{(p^{e}-1)(p^{f}-1)}$. Let $\gamma$ be a primitive element of $\F_q$. We compute the Gauss periods $\psi(\gamma^a C_0^{(N,q)})$, $0\leq a\leq N-1$, as follows.

\begin{align*}
\psi(\gamma^aC_0^{(N,q)})=\frac{1}{p^\ell-1}\sum_{x\in\F_{p^{e}}^\ast}\sum_{y\in\F_{p^{f}}^\ast}\xi_p^{\tr_{p^f}(y
\tr_{p^{3f}/p^f}(x\gamma^a))}=\frac{1}{p^\ell-1}\sum_{x\in\F_{p^{e}}^\ast}(p^{f}\delta_{\tr_{p^{3f}/p^{f}}(x\gamma^a)}-1),
\end{align*}
where
\[
\delta_{\tr_{p^{3f}/p^{f}}(x\gamma^a)}=\left\{
\begin{array}{ll}
1,&  \mbox{if $\tr_{p^{3f}/p^f}(x\gamma^a)=0$;}\\
0,&  \mbox{otherwise.}
 \end{array}
\right.
\]
Define
\[
W_a:=\{x\in\F_{p^{e}}\,|\,\tr_{p^{3f}/p^f}(x\gamma^a)=0\},
\]
and set $s_a=|W_a|$. Then we have
\[
\psi(\gamma^a C_0^{(N,q)})=
\frac{p^{f}(s_a-1)-(p^e-1)}{p^\ell-1}=\frac{p^{f}s_a-p^f-p^e+1}{p^\ell-1}.
\]
Since $W_a$ is an $\F_{p^\ell}$-subspace of $\F_{p^e}$, we have
$s_a=1,p^\ell,p^{2\ell},p^{3\ell}=p^e$.
Since a basis of $\F_{p^{e}}$ over $\F_{p^\ell}$ is also a basis of $\F_{p^{3f}}$ over $\F_{p^{f}}$, and $\gamma^a\neq 0$, it is impossible to have $W_a=\F_{p^{e}}$. Therefore
the Gauss periods $\psi(\gamma^a C_0^{(N,q)})$, $0\leq a\leq N-1$, take exactly three values
\[
\alpha_1=\frac{1-p^e}{p^\ell-1},\, \alpha_2=p^f+\frac{1-p^e}{p^\ell-1},\, \alpha_3=p^f(p^{\ell}+1)+\frac{1-p^e}{p^\ell-1}.
\]
By Lemma~\ref{rem_sz}, it is routine to compute that
\[
|I_1|=\frac{p^{3 \ell} + p^{2 f} - p^{2\ell + f} - p^{\ell + f}}{1 + p^\ell + p^{ 2 \ell}},\, |I_2|=p^{f}-p^{\ell},\, |I_3|=1.
\]
Since $|I_3|=1$, by Theorem~\ref{selfasso}, the subsets $\bigcup_{i\in I_j}C_i^{(N,q)}$, $j=1,2,3,$ give a three-class association scheme.

\subsection{Examples from index $2$ Gauss sums}\label{ind2} Let $q=p^f$, where $p$ is a prime and $f$ a positive integer. Let $N>1$ be a divisor of $q-1$. We now focus on the index $2$ case, that is, $[\Z_N^\ast:\langle p\rangle]=2$, or equivalently, $\ord_{N}(p)=\phi(N)/2$, where $\phi$ is Euler's phi function. In this case, the Gauss sums $G_q(\chi)$, where $\chi$ has order $N$, have been evaluated (cf. \cite{YX10}). In \cite{FX}, the authors used these Gauss sums to construct several new families of strongly regular graphs. In particular, they evaluated the Gauss periods in the index 2 case. The following theorem is a specialized version of Theorem 4.1 and Theorem 5.1 from \cite{FX}.

\begin{theorem} \label{ind2Thm}
\begin{itemize}
\item[(i)] {\rm (\cite[Theorem 4.1]{FX})}
Let $N=p_1\equiv 3\,(\mod{4})$ be a prime with $p_1>3$,
and let $p$ be a prime such that $\gcd{(p,N)}=1$ and $\ord_{N}(p)=(N-1)/2$.
Let $q=p^f$, where $f=(p_1-1)/2$. Then the Gauss periods $\psi(\gamma^a C_0^{(N,q)})$,
$a=0,1,\ldots,N-1$, take at most three values
\begin{equation}\label{eigenind2}
\alpha_1=\frac{-2 + p^{\frac{f-h}{2}}b(p_1-1)}{2p_1}, \,
\alpha_2=\frac{-2 + p^{\frac{f-h}{2}}cp_1-p^{\frac{f-h}{2}}b}{2p_1},
\, \alpha_3=\frac{-2 - p^{\frac{f-h}{2}}cp_1-p^{\frac{f-h}{2}}b}{2p_1},
\end{equation}
where $h$ is the class number of $\Q(\sqrt{-p_1})$, and $b$ and $c$ are integers
determined by $b,c\not\equiv 0\pmod{p}$, $4p^{h}=b^2+p_1c^2$,  and $bp^{\frac{f-h}{2}}\equiv -2\pmod{p_1}$.
\item[(ii)] {\rm (\cite[Theorem 5.1]{FX})}
Let $N=p_1p_2$, where $p_1$ and $p_2$ such that $p_1\equiv 1 (\mod{4})$ and $p_2\equiv 3 (\mod{4})$. Let $p$ be a prime such that $\ord_{p_1}(p)=p_1-1$, $\ord_{p_2}(p)=p_2-1$, $\ord_{p_1p_2}(p)=(p_1-1)(p_2-1)/2$.
Let $q=p^f$, where $f=(p_1-1)(p_2-1)/2$. Then the Gauss periods $\psi(\gamma^a C_0^{(N,q)})$,
$a=0,1,\ldots,N-1$,  take at most five values
\begin{align*}
&\alpha_1=\frac{-1 + \frac{1}{2} p^{\frac{f - h}{2}} (b + c p_1 p_2)}{N}, \,
\alpha_2=\frac{-1 + p^\frac{f}{2} (-\frac{1}{2} b p^\frac{-h}{2} (-1 + p_1) + p_1)}{N},\\
&\alpha_3=\frac{-1 + \frac{1}{2} p^{\frac{f - h}{2}} (b - c p_1 p_2)}{N}, \,
\alpha_4=\frac{-1 + p^\frac{f}{2} (-\frac{1}{2} b p^\frac{-h}{2} (-1 + p_2) - p_2)}{N},\\
&\alpha_5=\frac{-1 + p^\frac{f}{2} (p_1 + \frac{1}{2} b p^\frac{-h}{2} (-1 + p_1) (-1 + p_2) - p_2)}{N},
\end{align*}
where $h$ is the class number of $\Q(\sqrt{-p_1p_2})$, and $b$ and $c$ are integers
determined by $b,c\not\equiv 0\pmod{p}$, $4p^{h}=b^2+p_1p_2c^2$,  and $bp^{\frac{f-h}{2}}\equiv 2\pmod{p_1p_2}$.
\end{itemize}
\end{theorem}
From this theorem, we immediately have the following proposition.
\begin{proposition}\label{index2self}
\begin{itemize}
\item[(i)] With assumptions and notation the same as in Theorem~\ref{ind2Thm} (i), the Gauss periods $\psi(\gamma^a C_0^{(N,q)})$,
$a=0,1,\ldots,N-1$, take exactly three values which form an arithmetic progression if and only if  $p_1+9=4p^{h}$ and $\pm 3p^{(f-h)/2}\equiv -2\,(\mod{p_1})$.
\item[(ii)] With assumptions and notation the same as in Theorem~\ref{ind2Thm} (ii), the Gauss periods $\psi(\gamma^a C_0^{(N,q)})$,
$a=0,1,\ldots,N-1$,  take at most three values if $4 p^\frac{h}{2}\equiv 0\pmod{p_1+p_2}$ and $2 p^\frac{f}{2} (p_1 - p_2)/(p_1 + p_2)\equiv 2\pmod{p_1p_2}$. In particular, they take exactly three values forming an arithmetic
progression if and only if  $p_1p_2+9=4p^{h}$ and $\pm 3p^{(f-h)/2}\equiv 2\,(\mod{p_1p_2})$.
\end{itemize}
\end{proposition}
\proof
(i)  First we remark that from the explicit computations of the Gauss periods $\psi(\gamma^a C_0^{(N,q)})$ in the proof of Theorem 4.1 in \cite{FX}, we know that  if $\alpha_1, \alpha_2$ and $\alpha_3$ are distinct, then the Gauss periods take exactly three values, and $\alpha_1$ is taken precisely once.

It is clear that  $\alpha_1,\alpha_2,\alpha_3$ form an arithmetic progression if and only if  $b=\pm 3c$. Since $b,c\not\equiv 0 \pmod p$, we have $b=\pm 3c$ if and only if $c\in \{-1,1\}$ and $b=\pm 3$. It follows that the Gauss periods take exactly three values in arithmetic progression if and only if $p_1+9=4p^h$ and $\pm3p^{\frac{f-h}{2}}\equiv -2\pmod{p_1}$.

(ii)   Assume that $4 p^\frac{h}{2}\equiv 0\pmod{p_1+p_2}$ and $2 p^\frac{f}{2} (p_1 - p_2)/(p_1 + p_2)\equiv 2\pmod{p_1p_2}$. We set
\[
b= \frac{2 p^\frac{h}{2} (p_1 - p_2)}{p_1 + p_2} \mbox{\, and \, } c=\pm\frac{4 p^\frac{h}{2}}{p_1 + p_2}.
\]
Both $b$ and $c$ are integers, and they satisfy $4p^{h}=b^2+p_1p_2c^2$ and $bp^{\frac {f-h}{2}}\equiv 2\pmod {p_1p_2}$. Note that the above $b,c$ are all the integer solutions to $4p^{h}=b^2+p_1p_2c^2$ and $bp^{\frac {f-h}{2}}\equiv 2\pmod {p_1p_2}$. If $b= \frac{2 p^\frac{h}{2} (p_1 - p_2)}{p_1 + p_2} $ and $c=\frac{4 p^\frac{h}{2}}{p_1 + p_2}$, then $\alpha_1=\alpha_2$ and $\alpha_3=\alpha_4$. On the other hand, if $b= \frac{2 p^\frac{h}{2} (p_1 - p_2)}{p_1 + p_2} $ and $c=-\frac{4 p^\frac{h}{2}}{p_1 + p_2}$, then $\alpha_1=\alpha_4$ and $\alpha_2=\alpha_3$. In both cases,  the Gauss periods $\psi(\gamma^a C_0^{(N,q)})$, $a=0,1,\ldots,N-1,$ take at most three values $\alpha_1, \alpha_3, \alpha_5$ (also, from the computations in the proof of Theorem 5.1 in \cite{FX}, $\alpha_5$ occurs precisely once); in particular, these $\alpha_1,\alpha_3,\alpha_5$ form an arithmetic progression if and only if $p_1-p_2=\pm 6$ (i.e., $b=\pm 3 c$).  Since $b,c\not\equiv 0\pmod{p}$, we have $b=\pm 3c$ if and only if $c\in \{-1,1\}$ and $b=\pm 3$. It follows that n this case the Gauss periods take three values in arithmetic progression if and only if  $p_1p_2+9=4p^h$ and $\pm 3p^{\frac{f-h}{2}}\equiv 2\pmod{p_1p_2}$.
\qed

\begin{example}\label{exam:ind2}
There are only five examples satisfying the index 2 condition, and $p_1+9=4p^h$ and $\pm 3p^{(f-h)/2}\equiv -2\,(\mod{p_1})$ stated in Proposition~\ref{index2self} (i) for $p_1\le 20000$:
\[
(p_1,p,h)=(11,5,1),(23,2,3),(43,13,1),(67,19,1),(163,43,1).
\]
There are only two examples satisfying the index 2 condition, and $p_1p_2+9=4p^{h}$ and $\pm 3p^{(f-h)/2}\equiv 2\,(\mod{p_1p_2})$ stated in Proposition~\ref{index2self} (ii) for $p_1p_2\le 20000$:
\[
(p_1,p_2,p,h)=(5,11,2,4),(17,11,7,2).
\]
These results are obtained by a computer search.
\end{example}

\begin{remark}{\rm Let $q$ be a power of a prime $p$, $\gamma$ be a primitive element of $\F_q$, and $N>1$ be a divisor of $q-1$. In the semi-primitive case, i.e., the case where  $-1\in \langle p\rangle \,(\mod{N})$, it is well known that the Gauss periods $\psi(\gamma^aC_0^{(N,q)})$, $0\leq a\leq N-1$, take exactly two values. Note that the condition  $-1\in \langle p\rangle \,(\mod{N})$ does not involve the extension degree of $\F_q$ over $\Z_p$. Therefore, for any $e>1$, the Gauss periods corresponding to the subgroup of index $N$ of $\F_{q^e}^\ast$ also take exactly two values. One is thus led to the following question: are there examples of $(q,N)$, where $N|(q-1)$ and $N>1$, such that the Gauss periods $\psi(\gamma^aC_0^{(N,q)})$, $0\leq a\leq N-1$, take exactly three values, and for any $e>1$, the Gauss periods corresponding to the subgroup of index $N$ of $\F_{q^e}^\ast$ also take exactly three values?  The index $2$ case with $N=p_1$ gives a positive answer to this question. The reason is given below. Note that since $\Tr_{q/p}(x)=\Tr_{q/p}(x^p)$ for any $x\in \F_q$, each index set $I_i$ is invariant under the multiplication by $p$; in the index 2 case, it follows that each $I_i$ is a union of $\{0\},\langle p\rangle,-\langle p\rangle$. It is clear that this conclusion holds, irrelevant of the extension degree of $\F_q$ over $\Z_p$. Therefore, in this case, if the Gauss periods $\psi(\gamma^aC_0^{(N,q)})$, $0\leq a\leq N-1$, take exactly three values, then for any $e>1$, the Gauss periods corresponding to the subgroup of index $N$ of $\F_{q^e}^\ast$ also take exactly three values. Here, we should remark that the index $2$ case sometimes gives two-valued Gauss periods; all such possibilities are determined under the generalized Riemann hypothesis in \cite{sw}. Except for those examples of two-valued Gauss periods determined in \cite{sw}, the index 2 case with $N=p_1$ provides a positive answer to the question above.}
\end{remark}

\subsection{Computer search}\label{compu}
We conducted a computer search for examples of three-valued Gauss periods with the following restrictions: $p<300$, $p^f<2^{25}$, $3<N<1001$, $(p-1)|k=\frac{p^f-1}{N}$. The output is listed in Table \ref{table_cs}. Note that in Tabel~\ref{table_cs} we have removed the known examples given in the four subsections above because otherwise the table would take too much space. The multiplicities of the Gauss periods are given by the exponents; for example, in the first row of Table 1, $-7^{10}$ means that the Gauss periods $\eta_a$, $0\leq a\leq 18$, take the value $-7$ ten times. The AP column indicates whether the Gauss periods are in arithmetic progression or not, with ``$\circ$" meaning YES and ``$\times$" meaning No. The AS column indicates whether the index sets $I_j$, $j=1,2,3$, yield a three-class association scheme or not.
\begin{table}[h]
\begin{center}
\footnotesize
$$
\begin{array}{|c|c|c|c|c|c|}
\hline
 p & f& N&  \textup{Gauss periods} & AP&AS\\
\hline
\hline
11&3&19&-7^{10},4^6, 15^3& \circ&\times  \\\hline
 7& 7& 29&-414, -71^{21},272^{7}&\circ&\circ\\\hline
29&3& 67 &-13^{43},16^{18},45^6 &\circ&\times\\\hline
 37& 3& 67 & -21^{39}, 16^{18},53^{10}&\circ&\times  \\\hline
23& 3& 79 & -7^{58}, 16^{18}, 39^3 &\circ&\times  \\\hline
2& 11& 89&  -9^{11}, -1^{56},  7^{22}&\circ&\circ  \\\hline
5& 6& 93&   -7^{70},18^{20}, 43^3&\circ&\times \\\hline
37& 3& 201&   -7^{166}, 30^{32}, 67^{3}&\circ&\times   \\\hline
67& 3& 217&   -21^{159}, 46^{48},113^{10}&\circ&\times  \\\hline
2& 18& 219&  -19^{163}, 45^{47}, 109^{9} &\circ&\times \\\hline
61& 3& 291&   -13^{235}, 48^{50}, 109^6&\circ&\times  \\\hline
79& 3& 301&    -21^{231},58^{60},137^{10} &\circ&\times  \\\hline
 83& 3& 367&    -19^{292},64^{66}, 147^9 &\circ&\times \\\hline
 11& 6& 399&   -37^{295},  84^{86},205^{18}&\circ&\times   \\\hline
\end{array}
\begin{array}{||c|c|c|c|c|c|}
\hline
 p & f& N&  \textup{Gauss periods} & AP&AS\\
\hline
\hline
53& 3& 409&  -7^{358}, 46^{48}, 99^3 &\circ &\times  \\\hline
 139& 3& 499&   -39^{378},  100^{102},239^{19}&\circ&\times   \\\hline
 137& 3& 511&   -37^{391},  100^{102},237^{18}&\circ&\times   \\\hline
 109& 3& 571&   -21^{471},  88^{90},197^{10}&\circ&\times   \\\hline
67& 3& 651&   -7^{586}, 60^{62}, 127^{3}&\circ&\times \\\hline
 11& 6& 703&   -21^{591},  100^{102},221^{10}&\circ&\times   \\\hline
 149& 3& 721&   -31^{586},  118^{120},267^{15}&\circ&\times   \\\hline
 11& 6& 777&   -19^{661},  102^{113},343^{3}&\times&\times   \\\hline
 5& 9& 829&   -19^{712},  106^{108},231^{9}&\circ&\times   \\\hline
 107& 3& 889&   -13^{787},  94^{96},201^{6}&\circ&\times   \\\hline
79& 3& 903&   -7^{826}, 72^{74},  151^{3}&\circ&\times \\\hline
 17& 6& 921&   -91^{676},  198^{200},487^{45}&\circ&\times   \\\hline
3& 12& 949&  -7^{870},  74^{76}, 155^{3}&\circ&\times  \\\hline
 113& 3& 991&   -13^{883},  100^{102},213^{6}&\circ&\times   \\\hline
\end{array}
$$
\caption{ Computer search results for $p<300$, $p^f<2^{25}$, $6<N<1001$, $N|\frac{p^f-1}{p-1}$ except for the known examples given in Subsections~\ref{quadric}-\ref{ind2}}.
\label{table_cs}
\end{center}
\end{table}

Furthermore, Corollary~\ref{co:self} makes it possible to search for $(p,f,N)$ such that the Gauss periods corresponding to the subgroup of index $N$ of $\F_q^\ast$, $q=p^f$, take exactly three values.


We will run the following algorithm to search for triples $(p,f,N)$ satisfying the conditions in Corollary~\ref{co:self}:  (i) $t(vs-ur)+1\equiv 0\,(\mod{N})$, (ii) $(N-1)q+
t^2(vs-ur)^2=(u^2r+v^2s)t^2N$, and (iii) $u=v=1$ and $r=1$ or $s=1$. Put $g=s-r$ and $h=r+s$. In this case, we have $h=|g|+2$. The
algorithm goes as follows:
\begin{itemize}
\item[(1)]
For any positive integers $N$ and $h$ with $1<h<N$, compute $(Nh-(h-2)^2)/(N-1)$ in order to know $q/t^2$.
\item[(2)] If this value is a prime power, say $p^w$, then compute
the order of $p$ modulo $N$, call it $f'$, and the largest positive integer $p^{\theta'}$ dividing $G_{p^{f'}}(\chi)$ for all nontrivial characters $\chi$ of
exponent $N$ of $\F_{p^{f'}}^\ast$.
\item[(3)] Check
whether $f'-2\theta'$ divides $w$. Set $d=w/(f'-2\theta')$ and $t=p^\theta=p^{d\theta'}$. Then,
check whether $(h-2)t+1\equiv 0\,(\mod{N})$ or
$-(h-2)t+1\equiv 0\,(\mod{N})$ holds.
\end{itemize}
We run the above algorithm for all $N<5000$ using a computer.  Note that $p$ is determined as the unique prime factor of $(Nh-(h-2)^2)/(N-1)$ in Steps (1) and  (2), and $f$ is determined as $f=df'$ in the steps (2) and (3). We find three quadruples (for convenience we  give the value of $\theta$ also) satisfying the conditions of Corollary~\ref{co:self}:
\begin{equation}\label{spo}
(p,f,N,\theta)=(7,7,29,3),(13,13,53,6),(2,36,247,15).
\end{equation}
By Theorem~\ref{selfasso}, we obtain three new self-dual three-class association schemes from the three quadruples above. These self-dual 3-class association schemes are different from the examples obtained in Subsections~\ref{primaryEX} and \ref{ind2}.


As a counterpart of Conjecture 4.4 in \cite{sw}, we have the following conjecture.

\begin{conj} Let $q$ be a power of a prime $p$, $\gamma$ be a primitive element of $\F_q$, and $N>1$ be a divisor of $q-1$.
The Gauss periods $\psi(\gamma^a C_0^{(N,q)})$, $a=0,1,\ldots,N-1$, take
exactly three rational values in arithmetic progression, and one of three values occurs exactly once, if and only if the Gauss periods arise from the examples in Subsection~\ref{primaryEX}, or from Example~\ref{exam:ind2}, or from one of the sporadic cases listed in (\ref{spo}).
\end{conj}

\section{Concluding remarks}
In this paper, we study the problem of when the Gauss periods take exactly three rational values. Also, we give constructions of related combinatorial structures such as circulant weighing matrices and association schemes.

We have found five infinite classes of three-valued  Gauss periods
listed in Table~\ref{Tab2}.
\begin{table}[h]
$$
\begin{array}{|c||c|c|c|c|}
\hline
\mbox{parameters}&\mbox{AP}&\mbox{AS}&\mbox{CW}&\mbox{ref} \\
\hline \hline
\mbox{$p=2$, $q=p^{6f}$, $N=\frac{p^{3f}-1}{p^f-1}$} &\circ&\circ&\circ&\mbox{Subsec.~\ref{quadric}}\\
\hline
\mbox{$p$ odd, $q=p^{6f}$, $N=\frac{p^{3f}-1}{p^f-1}$}
&\circ&\times&\circ&\mbox{Subsec.~\ref{quadric}}\\
\hline
\mbox{$q=p^{3f}$,$N=\frac{p^{3f}-1}{p^f-1}$, $\ord_{3(p^f-1)}(p^f)=3$}
&\circ&\circ&\times&\mbox{Subsec.~\ref{primaryEX}}\\
\hline
\mbox{$q=p^{fe}$, $\frac{p^{fe}-1}{N}\,|\,p^f-1$, \, $\frac{(p^f-1)N}{p^{fe}-1}\,|\,\frac{p^f-1}{p-1}$, $Cay(\F_q,C_0^{(\frac{(p^f-1)N}{p^{fe}-1},p^f)})$ is an SRG}
&\star&\circ&\times&\mbox{Subsec.~\ref{twoGauss}}\\
\hline
\mbox{$q=p^{\lcm{(e,f)}}=p^{3f}$, $e/\gcd{(e,f)}=3$,  $C_0^{(N,q)}=\F_{p^e}^\ast\cdot \F_{p^f}^\ast$}
&\times&\circ&\times&\mbox{Subsec.~\ref{sec:subf}}\\
\hline
\mbox{$N=p_1$, $[\Z_N^\ast,\langle p\rangle ]=2$, $f=e(N-1)/2$ for any $e\in \N$}& \star&\circ&\times&\mbox{Subsec.~\ref{ind2}}\\
\hline
\mbox{$N=p_1p_2$, $[\Z_N^\ast:\langle p\rangle ]=2$, $f=\phi(N)/2$}& \star&\circ&\times&\mbox{Subsec.~\ref{ind2}}\\
\hline
\end{array}
$$
\caption{\label{Tab2}Known examples of three-valued Gauss periods
}
\end{table}
(The meaning of ``AP,'' ``AS'' are the same as in Table 1.  Here ``CW'' indicates whether $I_1-I_3$ gives a circulant weighing matrix or not. The symbols ``$\star$'' means that the class includes some examples satisfying the condition.) Furthermore, we obtained several sporadic examples of
three-valued Gauss periods as given in Subsection~\ref{compu}.

We conclude the paper by listing some problems for future work.
\begin{itemize}
\item Classify all $(p,f,N)$ which lead to three-valued Gauss periods.  A less challenging task is to find other infinite classes  of
three-valued Gauss periods not listed in Table~\ref{Tab2}.
\item Determine when three-valued Gauss periods take three values in arithmetic progression. (Then, by Proposition~\ref{circu} one will be able to characterize when $I_1-I_3$ forms a circulant weighing matrix.)
\item Determine when the index sets $I_1,I_2,I_3$ yield a three-class association scheme if the Gauss periods take exactly three values.
\end{itemize}

\section*{Acknowledgement}
The authors would like to thank both reviewers for their comments and constructive suggestions. In particular, we thank one of the reviewers who gave a short proof of Theorem~\ref{selfasso}, which is the proof presented here in this paper.

\end{document}